\documentstyle[12pt]{article}
\def\Hom{\mathop{\rm Hom}\nolimits}
\def\End{\mathop{\rm End}\nolimits}
\def\Ker{\mathop{\rm Ker}\nolimits}
\def\Im{\mathop{\rm Im}\nolimits}
\def\tr{\mathop{\rm tr}\nolimits}
\def\e{{\varepsilon}}
\def\t{{\tau}}
\def\th{{\theta}}
\def\b{{\beta}}
\def\la{{\langle}}
\def\ra{{\rangle}}
\def\v{{\vert}}
\def\V{{\Vert}}
\def\arr{{\longrightarrow}}
\def\a{{\alpha}}

\def\l{{\lambda}}
\def\i{{\infty}}
\def\q{{$\quad\bullet$}}

\def\A{{\cal A}}
\def\H{{\cal H}}
\def\K{{\cal K}}
\def\B{{\cal B}}

\def\C{{$C^*$-algebra}}

\date{22 January 1995}
\author{V.~M.~Manuilov}
\title{Diagonalization of compact operators in Hilbert
       modules over  $C^*$-algebras of real rank zero}
\begin{document}

\maketitle
\begin{abstract}
It is known that the classical Hilbert--Schmidt theorem can be generalized
to the case of compact operators in Hilbert $\A$-modules $\H_\A^*$ over
a $W^*$-algebra of finite type, i.e. compact operators in $\H_\A^*$ under
slight restrictions can be diagonalized over $\A$. We show that if $\B$ is
a weakly dense $C^*$-subalgebra of real rank zero in $\A$ with some
additional property then the natural
extension of a compact operator from $\H_\B$ to $\H_\A^*\supset \H_\B$ can
be diagonalized with diagonal entries being from the $C^*$-algebra $\B$.
\end{abstract}

\section{Introduction}
Let $\A$ be a \C. We consider Hilbert $\A$-modules over $\A$ \cite{pa1},
i.e. (right)
$\A$-modules $\cal M$ together with an $\A$-valued inner product $\la\cdot,
\cdot\ra : {\cal M}\times{\cal M}\arr \A$ satisfying the following
conditions:
\begin{enumerate}
\item
$\la x,x\ra\geq 0$ for every $x\in {\cal M}$ and $\la x,x\ra =0$ iff $x=0$,
\item
$\la x,y\ra =\la y,x\ra^*$ for every $x,y\in {\cal M}$,
\item
$\la\cdot,\cdot\ra$ is $\A$-linear in the second argument,
\item
$\cal M$ is complete with respect to the norm $\V x\V^2=\V\la x,x\ra\V_\A$.
\end{enumerate}
By $\cal M^*=\Hom_\A(\cal M;\A)$ we denote the $\A$-module dual to $\cal M$.
Let $\H_\A$ be a right Hilbert $\A$-module
of sequences
$a=(a_k)$, $a_k\in \A$, $k\in {\bf N}$ such that the series $\sum a^*_ka_k$
converges in $\A$ in norm with the standard basis $\{e_k\}$ and
let $L_n(\A)\subset \H_\A$ be a submodule generated by the elements
$e_1,\ldots,e_n$ of the basis. An inner $\A$-valued product on
module $\H_\A$ is given by $\la x,y\ra=\sum x^*_ky_k$ for $x,y\in\A$.
A bounded operator $\K:\H_\A\arr \H_\A$ is
called compact~\cite{ka}~\cite{mf}, if it possesses an adjoint operator and
lies in the norm closure of the
linear span of operators of the form $\th_{x,y}$, $\th_{x,y}(z)=x\la
y,z\ra$, $x,y,z\in \H_\A$. From now on we  suppose that the compact
operator $\K$ is strictly positive, i.e. operator $\la \K x,x\ra$ is
positive in $\A$ and $\Ker \K=0$.
It is known~\cite{pa2} that in the case when $\A$ is a $W^*$-algebra the
inner product can be naturally prolonged to the dual module $\H^*_\A$.

\medskip
{\bf Definition 1.1.} \ Let $\A$ be a $W^*$-algebra. We call an operator
$\K$ {\it diagonalizable\/} if
there exist a set
$\{x_i\}$ of
elements in $\H^*_\A$ and a set of operators $\l\in\A$ such that
\begin{enumerate}
\item
$\{x_i\}$ is orthonormal, $\la x_i,x_j\ra=\delta_{ij}$,
\item
$\H_\A^*$ coincides with the $\A$-module $\cal M^*$ dual to the module
$\cal M$ generated by the set $\{x_i\}$,
\item
$\K x_i=x_i\l_i$,
\item
for any unitaries $u_i,\,u_{i+1}\in\A$ we have an operator inequality
\begin{equation}\label{ord}
u_i^*\l_iu_i\ge u^*_{i+1}\l_{i+1}u_{i+1}.
\end{equation}
\end{enumerate}

We call the elements $x_i$ ``{\it eigenvectors\/}'' and the operators $\l_i$
``{\it eigenvalues\/}'' for the operator $\K$. It must be noticed that the
``eigenvectors'' and ``eigenvalues'' are defined not uniquely.

\smallskip
The problem of diagonalizing operators in Hilbert modules was initiated by
R.~V.~Kadison in~\cite{kad} and was studied in
different settings in~\cite{gp},\cite{mur},\cite{fm},\cite{st} etc.
In~\cite{man1},\cite{man2} we have proved the following

\medskip
{\bf Theorem 1.2.}\ {\it If $\A$ is a finite $\sigma$-finite $W^*$-algebra
then a compact strictly positive operator $\K$ can be diagonalized and its
``eigenvalues'' are defined uniquely up to unitary equivalence.}

\medskip
It is well known that in the commutative case, i.e. for ${\cal C}=C(X)$ being
a commutative \C, compact operators cannot be diagonalized inside $\H_{\cal
C}$ but it becomes possible if we pass to a bigger module over a bigger
$W^*$-algebra $L^\infty(X)\supset\cal C$. It leads us to the following

\medskip
{\bf Definition 1.3.}\ Let $\cal C$ be a
\C~ admitting a weakly dense inclusion in a finite $\sigma$-finite
$W^*$-algebra $\A$ and let $\K$
be a compact strictly positive operator in $\H_{\cal C}$. We can naturally
extend $\K$ to the bigger module $H^*_\A$ where it will remain compact and
strictly positive and by the theorem 1.2 it can be
diagonalized in this module. We call a \C~ $\cal C$ admitting {\it weak
diagonalization} if the diagonal entries for any $\K$ in $\H^*_\A$ can be
taken from $\cal C$ instead of $\A$.

\medskip
{\bf Problem.} Describe the class of \C s admitting weak diagonalization.

\medskip
Throughout this paper we denote by $\A$ a finite $\sigma$-finite
$W^*$-algebra. Denote by ${\cal Z}=C(Z)$ the center of $\A$ and by
$T$ the standard exact center-valued trace defined on $\A$,
$T({\bf 1})=1$. Suppose that for a $C^*$-subalgebra $\B$ of $\A$ the
following condition holds:
\begin{itemize}
\item[({\raisebox{-3pt}{\bf *}})] for any two projections $p,q\in\B$ there
exist in $\B$ equivalent (in $\B$) projections $r_p\sim r_q$,
$r_p\leq p$, $r_q\leq q$ such that
$T(r_p)=T(r_q)=\min\{T(p)(z),T(q)(z)\}$,
$z\in Z$.
\end{itemize}
The purpose of this paper is to show that the class of \C s admitting
weak diagonalization contains
real rank zero
weakly dense $C^*$-subalgebras of finite $\sigma$-finite
$W^*$-algebras with the property ({\raisebox{-3pt}{\bf *}}).
Recall that real rank zero ($RR(\B)=0$) means~\cite{bp} that every
selfadjoint operator in $\B$ can be approximated by operators with finite
spectrum, i.e. having the form $\sum\a_ip_i$, where
$p_i\in\B$ are selfadjoint mutually orthogonal projections and
$\a_i\in{\bf R}$. By~\cite{bp} we have in this case also
$RR(\End_\B(L_n(\B)))=0$.

\section{Continuity of ``eigenvalues''}
\setcounter{equation}{0}
For the further we need to establish some continuity properties of the
``eigenvalues'' of compact operators in modules over $W^*$-algebras.

\bigskip
{\bf Lemma 2.1.}\ {\it Let $\K_1=\sum\a_l^{(1)}P^{(1)}_l$, $\K_2=\sum\a_l^{(2)}
P^{(2)}_l$ be
strictly positive operators in $L_n(\A)$ with finite spectrum and let
$\V \K_1-\K_2\V<\e$. Then
\begin{enumerate}
\item
one can find a unitary $U$ in $L_n(\A)$ such that it maps the
``eigenvectors'' of $\K_2$ to the ``eigenvectors'' of $\K_1$ and
$\V U^*\K_1U-\K_2\V<\e$,
\item
``eigenvalues'' $\{\l_i^{(r)}\}$ of operators $\K_r$ $(r=1,\,2)$ can be
chosen in such a way that $\V\l_i^{(1)}-\l_i^{(2)}\V<\e$.
\end{enumerate}}

\medskip
{\bf Proof.}\ As the algebra $\A$ can be decomposed into a direct integral
of finite factors, so it is sufficient to prove the lemma for the case when
$\A$ is a type ${\rm II}_1$ factor (for type ${\rm I}_n$ factors lemma is
trivial). Denote by $E_\K(\l)$ the spectral
projection
for the operator $\K$ corresponding to the set $(-\i,\l)$. If $\t$ is an
exact finite trace on $\A$, it can be prolonged to the (infinite) trace
$\bar\t=\tr\otimes\t$ on the algebra $\End_\A(\H_\A^*)$ and to the finite
trace on a lesser algebra $\End_\A(L_n(\A))$ where we have $\bar\t({\bf 1})=
n$. Put $$
\e_\K(\a)=\inf_{\bar\t(E_\K(\l))\ge\a}\l, \quad 0\le\a\le n. $$
As it is shown in~\cite{mn} (the continuous minimax principle) one has
\begin{equation}\label{eps}
\e_\K(\a)=\inf_{P\in{\cal P},
\,\bar\t(P)\ge\a}\Bigl\lbrace\sup_{\xi\in\Im P,\,
\V\xi\V=1}(\K\xi,\xi)\Bigr\rbrace,
\end{equation}
where $(\cdot,\cdot)$ denotes an inner product in a Hilbert space where the
algebra $\End_\A(L_n(\A))$ is represented and ${\cal P}$ denotes the set of
projections in $\End_\A(L_n(\A))$.
It follows from (\ref{eps}) that if $\V \K_1-\K_2\V<\e$, then
\begin{equation}\label{eps2}
\v\e_{\K_1}(\a)-\e_{\K_2}(\a)\v<\e.
\end{equation}
Let $Q_i^{(r)}$ be projections on the ``eigenvectors'' $x_i^{(r)}$ of the
operators
$\K_r$, corresponding to the maximal ``eigenvalues'' $\l^{(r)}_i$,
$\bar\t(Q_i^{(r)})=1$. For two divisions $\{P_l^{(1)}, Q_i^{(1)}\}$ and
$\{P_l^{(2)}, Q_i^{(2)}\}$
of unity given by decompositions of $\K_1$ and $\K_2$ we can construct a
finer
division of unity. By~\cite{tak} there exist sets of mutually orthogonal
projections $R^{(r)}_m\in\End_\A(L_n(\A))$ such that
\begin{enumerate}
\item
$\bigoplus_mR^{(r)}_m=1$,
\item
$\bar\t(R^{(1)}_m)=\bar\t(R^{(2)}_m)$,
\item
for every $m$ we have $R^{(r)}_m\le Q_i^{(r)} \ {\rm or}\ R^{(r)}_m\le
P^{(r)}_j$ for some $i$ or $j$.
\end{enumerate}
Then (after renumbering) one can write the operators $\K_r$ in the form
$\K_r=\sum\a^{(r)}_mR^{(r)}_m$ with $\a^{(r)}_1\le\a^{(r)}_2\le\ldots$,
$\a^{(r)}_m\in\bf R$.
It makes possible to
define a unitary $U:L_n(\A)\arr L_n(\A)$ such that
\begin{equation}\label{U}
U(\Im R^{(2)}_m)=\Im R^{(1)}_m,
\end{equation}
hence $U(\Im Q_i^{(2)})=\Im Q_i^{(1)}$ so $U$ maps the
$\A$-modules generated by the ``eigenvectors'' $x_i^{(2)}$ into the modules
generated by $x_i^{(1)}$, hence $Ux_i^{(2)}=x_i^{(1)}\cdot
u_i=\bar x_i^{(1)}$ for some unitaries $u_i\in\A$.
Put
$$n(\a)=
\min \{n\v\bar\t(\bigoplus_{m\ge n}R^{(r)}_m)\ge\a\}.$$ Then
$\e_{\K_r(\a)}=\a^{(r)}_{n(\a)}$ and it follows from (\ref{eps2}) that
$\v\a^{(1)}_{n(\a)}-\a^{(2)}_{n(\a)}\v<\e$.
But changing $\a$ we obtain that
\begin{equation}\label{eps3}
\v\a^{(1)}_m-\a^{(2)}_m\v<\e
\end{equation}
for all $m$. Taking $\a=1$ (then $i=1$) we have $$\K_r\v_{\Im
Q_1^{(r)}}=\Lambda^{(r)}_1=\sum_{m\ge n(1)}\a^{(r)}_mP^{(r)}_m.$$ From
(\ref{U}) and (\ref{eps3})  we conclude that
\begin{equation}\label{a}
\V U^*\Lambda^{(1)}_1U-\Lambda^{(2)}_1\V=\V\sum_{m\ge n(1)}
(\a^{(1)}_m-\a^{(2)}_m)P^{(2)}_m\V
\le\e\V\bigoplus_{m\ge n(1)}P^{(2)}_m\V=\e.
\end{equation}
Choosing appropriate $\l_1^{(r)}$ to satisfy the conditions
$\Lambda^{(1)}_1\bar x^{(1)}_1=
\bar x^{(1)}_1\l^{(1)}_1$ and $\Lambda^{(2)}_1x^{(2)}_1=x^{(2)}_1\l^{(2)}_1$
we obtain the estimate
\begin{equation}\label{b}
\V\l^{(1)}_1-\l^{(2)}_1\V<\e.
\end{equation}
By the same way  estimates (\ref{a}), (\ref{b}) can be
obtained for all $i$ and it proves the lemma.\q

\medskip
{\bf Corollary 2.2.}\ {\it Let $\K_r:\H_\A\arr \H_\A$, $r=1,\,2$ be compact
strictly positive operators and let
$\V \K_1-\K_2\V<\e$. Then
\begin{enumerate}
\item
one can find a unitary $U$ in $\H_\A^*$ such that it maps the
``eigenvectors'' of $\K_2$ to the ``eigenvectors'' of $\K_1$ and
$\V U^*\K_1U-\K_2\V<\e$,
\item
``eigenvalues'' $\{\l_i^{(r)}\}$ of operators $\K_r$ $(r=1,\,2)$ can be
chosen in such a way that $\V\l_i^{(1)}-\l_i^{(2)}\V<\e$.
\end{enumerate}}

\medskip
{\bf Proof.}\ Let $L^{(r)}_n(\A)\in \H_\A^*$ denotes the Hilbert submodule
generated by the first $n$ ``eigenvectors'' of the operator $\K_r$,
$L^{(r)}_n(\A)\cong L_n(\A)$.
It was shown in~\cite{man2} that the orthogonal complement to such
submodule is isomorphic to $\H_\A^*$ and
the norm of restriction of compact operator $\K_r$ on the
orthogonal complement to $L^{(r)}_n(\A)$ in $\H_\A^*$ tends to zero,
henceforth it is
sufficient to consider only the case of operators in $L_n(\A)$ and there
one can approximate these operators by operators with finite spectrum.\q

\section{Case of $RR(\B)=0$}
\setcounter{equation}{0}
In this section we show that \C s of real rank zero with the property
({\raisebox{-3pt}{\bf *}}) admit weak diagonalization.

\bigskip
{\bf Theorem 3.1.}\ {\it
Let $\B$ be a weakly dense $C^*$-subalgebra in $\A$ with the property
({\raisebox{-3pt}{\bf *}}) and let $RR(\B)=0$.
If $\K$ is a compact strictly positive operator in the $\B$-module
$\H_\B$ then the ``eigenvalues'' $\{\l_i\}$ of diagonalization of
the natural prolongation of $\K$ to the $\A$-module $\H^*_\A$ can be chosen
in a way that $\l_i\in\B$ would hold.}

\medskip
{\bf Proof} is based on the results of S.~Zhang~\cite{zh}. By~\cite{bp},%
\cite{zh}
the operator $\K$ can be approximated by operators $\K_n\in\End_\B(L_n(\B))$
with finite spectrum. By~\cite{zh}, corollary 3.5 there exist such
unitaries $U_n\in\End_\B(L_n(\B))$ that the operators
$$U^*_n\K_nU_n=\left(\begin{array}{ccc}
          \l^{(n)}_1 & &0\\
            &\ddots& \\
          0 & &\l^{(n)}_n
         \end{array}\right)$$
are diagonal and $\l_i^{(n)}\in\B$ are operators with finite spectrum.
Show that due to the property ({\raisebox{-3pt}{\bf *}}) by an appropriate
choice of such $U_n$ one can make the condition
(\ref{ord}) valid for ``eigenvalues'' $\{\l_i^{(n)}\}$.
Let $\l_a=\sum\a_kq_k$,
$\l_b=\sum\b_lr_l$ where $q_k,\, r_l\in\B$ are projections and suppose that
$a<b$ but for some $m$ and $n$ inequality $\b_m>\a_n$ holds. Using the
possibility to diagonalize projections~\cite{zh} we can find projections
$s_l\in\B$ equivalent to $r_l$ and such that $s_l=\oplus_ks^{(l)}_k$ and
$s^{(l)}_k\le q_k$. Then put
$$\l'_a=\sum_{k\not= n}\a_kq_k\oplus\sum_{l\not=
m}\a_ns^{(l)}_n\oplus\b_ns^{(m)}_n,$$
$$\l'_b=\sum_{l\not= m}\b_ls_l\oplus\sum_{k\not=
n}\b_ks^{(m)}_k\oplus\a_ns^{(m)}_n$$
and notice that the operators
{\small       $\left(\begin{array}{ccc}
             \l_1&& 0\\
            &\ddots& \\
            0& & \l_n      \end{array}\right)$ }
and
{\small       $\left(\begin{array}{ccc}
             \l'_1&& 0\\
            & \ddots& \\
             0& & \l'_n      \end{array}\right)$ }
are unitarily equivalent. After repeating this procedure for all cases when
$\b_l>\a_k$ we obtain validity of (\ref{ord}) for $\l'_a$ and $\l'_b$.
By the same
way we can order all ``eigenvalues'' of $\K_n$ remaining in $\B$. But by the
property (\raisebox{-3pt}{\bf *}) if $\V \K_n - \K_{n-1}\V<\e_n$ then one
can find such unitaries $u_{i,n}$ in $\B$ that
\begin{equation}\label{en}
\V u_{i,n}^*\l^{(n)}_iu_{i,n}-\l_i^{(n-1)}\V<\e_n
\end{equation}.
Then $u_{i,n}^*\l^{(n)}_iu_{i,n}\in\B$.
Taking a
subsequence of $\{\K_n\}$ if necessary we can take in (\ref{en})
$\e_n=\frac{1}{2^n}$. Then the sequence $$\bar\l^{(1)}_i=\l^{(1)}_i,\
\bar\l^{(2)}_i=
u_{i,2}^*\l^{(2)}_iu_{i,2},\
\bar\l^{(3)}_i=u_{i,3}^*u_{i,2}^*\l^{(3)}_iu_{i,2}u_{i,3},\,\ldots$$
is fundamental
in $\B$. Denote its limit by $\bar\l_i\in\B$.
By the corollary 2.2 for all $\K_n$ we can find unitaries $U_n$ which
map the first $n$ ``eigenvectors'' of $\K$ to ``eigenvectors'' of $\K_n$.
Put $\K'_n=U^*_n\K_nU_n\in\End_\A(\H_\A^*)$. Then we have
\begin{equation}\label{n}
\K'_nx_i=x_i\bar\l^{(n)}_i
\end{equation}
and $\V \K'_n-\K\V\to 0$. Taking limit in (\ref{n}) we obtain
$\K x_i=x_i\bar\l_i$, hence $\bar\l_i$ are ``eigenvalues'' of $\K$.\q

\bigskip
Notice that the condition (\raisebox{-3pt}{\bf *}) is necessary for a
$C^*$-algebra to have the weak diagonalization property. Indeed if
$\K$ is a direct sum of two projections, {\small $\K=\left(\begin{array}{cc}
p&0\\0&q\end{array}\right)$} then the ``eigenvalues'' of $\K$ can be
ordered only if the ``common part'' of $1-p$ and $q$ lies in $\B$.

\bigskip
{\bf Remark.}\
%
In the case of  \C s $A_\th$ of irrational rotation one has
$RR(A_\th)=0$ (cf~\cite{ce}) and the property ({\raisebox{-3pt}{\bf *}})
is valid, so the theorem 3.1 gives the
answer to the problem of~\cite{man2} where we have considered the
Schr\"odinger operator in magnetic field with irrational magnetic flow.
It is known that this operator can be viewed
as an operator acting
in a Hilbert $A_\th$-module. As we can imbed $A_\th$ in a type ${\rm II}_1$
factor $\A$ as a weakly dense subalgebra~\cite{br} so we can diagonalize this
operator
in a Hilbert $\A$-module. The present paper shows that the ``eigenvalues''
of this operator can be chosen to be elements of $A_\th$. So this situation
is a
noncommutative analogue of the case $\th=1$ when the corresponding operator
can be diagonalized over $W^*$-algebra $L^\i({\bf T}^2)$ but the diagonal
elements lie in a lesser \C~$C({\bf T}^2)$. Notice that in case of rational
$\th$ this operator is also diagonalizable.

\bigskip
{\bf Acknowledgement.}\
This work was partially
supported by the Russian Foundation for Fundamental Research (grant
\mbox{N 94-01-00108-a)} and the International Science Foundation
(grant N MGM000).
I am indebted to M.~Frank, A.~A.~Irmatov, A.~S.~Mishchenko
and E.~V.~Troitsky for helpful discussions.

\vspace{2.5cm}
\noindent
V.~M.~Manuilov \\*
Moscow State University \\*
Moscow, 119899, Russia  \\*
E-mail:manuilov@mech.math.msu.su

\end{document}